\documentclass[12pt,oneside,reqno]{amsart}
\usepackage{txfonts}
\usepackage{bbm}
\usepackage{amsmath}
\usepackage{graphicx}
\usepackage{mathrsfs}
\usepackage{stmaryrd}
\usepackage{amsfonts}
\usepackage{enumerate,amsmath,amssymb,amsthm}

\numberwithin{equation}{section}

\newcommand{\be}{\begin{eqnarray}}
\newcommand{\ee}{\end{eqnarray}}
\newcommand{\ce}{\begin{eqnarray*}}
\newcommand{\de}{\end{eqnarray*}}
\newtheorem{theorem}{Theorem}[section]
\newtheorem{lemma}[theorem]{Lemma}
\newtheorem{remark}[theorem]{Remark}
\newtheorem{definition}[theorem]{Definition}
\newtheorem{proposition}[theorem]{Proposition}
\newtheorem{Examples}[theorem]{Example}
\newtheorem{corollary}[theorem]{Corollary}

\def\eps{\varepsilon}

\def\[{{\Big[}}
\def\]{{\Big]}}
\def\<{{\langle}}
\def\>{{\rangle}}
\def\({{\Big(}}
\def\){{\Big)}}

\def\bx{{\mathbf{x}}}

\def\dif{{\mathord{{\rm d}}}}

\def\no{\nonumber}
\def\={&\!\!=\!\!&}
\def\bt{\begin{theorem}}
\def\et{\end{theorem}}
\def\bl{\begin{lemma}}
\def\el{\end{lemma}}
\def\br{\begin{remark}}
\def\er{\end{remark}}

\def\bd{\begin{definition}}
\def\ed{\end{definition}}
\def\bp{\begin{proposition}}
\def\ep{\end{proposition}}
\def\bc{\begin{corollary}}
\def\ec{\end{corollary}}
\def\bx{\begin{Examples}}
\def\ex{\end{Examples}}

\def\cN{{\mathcal N}}

\def\cS{{\mathcal S}}

\def\mD{{\mathbb D}}
\def\mE{{\mathbb E}}

\def\mG{{\mathbb G}}
\def\mH{{\mathbb H}}

\def\mN{{\mathbb N}}

\def\mR{{\mathbb R}}

\def\mT{{\mathbb T}}

\def\bP{{\mathbf P}}

\def\sA{{\mathscr A}}
\def\sB{{\mathscr B}}

\def\sF{{\mathscr F}}

\def\sL{{\mathscr L}}

\def\geq{\geqslant}
\def\leq{\leqslant}

\def\div{\mathord{{\rm div}}}

\def\bP{{\mathbf P}}

\allowdisplaybreaks

\begin{document}

\title{Invariance measures of stochastic 2D Navier-Stokes equations driven by $\alpha$-stable processes}

\date{}
\author{Zhao Dong, Lihu Xu, Xicheng Zhang}

\address{Zhao Dong: Institute of Applied Mathematics, 
Academy of Mathematics and Systems Sciences, Academia Sinica, P.R.China\\ 
Email: dzhao@amt.ac.cn}

\address{Lihu Xu: TU Berlin, Fakult\"at II, Institut f\"ur Mathematik, 
Str$\alpha\beta$e Des 17. Juni 136, D-10623 Berlin, Germany\\
Email: xu@math.tu-berlin.de}

\address{Xicheng Zhang:
School of Mathematics and Statistics, Wuhan University,
Wuhan, Hubei 430072, P.R.China\\
Email: XichengZhang@gmail.com
 }

\begin{abstract}
In this note we prove the well-posedness for stochastic 2D Navier-Stokes equation driven by general L\'evy processes
(in particular, $\alpha$-stable processes), and obtain the existence of invariant measures.
\end{abstract}

\maketitle
\rm

\section{Introduction and Main Result}

In this article we are concerned with the following stochastic 2D Navier-Stokes equation in torus $\mT^2=(0,1]^2$:
\begin{align}
\dif u_t=[\Delta u_t-(u_t\cdot\nabla)u_t+\nabla p_t]\dif t+\dif L_t,
\ \ \div u_t=0,\ \ u_0=\varphi\in \mH^0,\label{NSE}
\end{align}
where $u_t(x)=(u^1_t(x),u^2_t(x))$ is the 2D-velocity field, $p$ is the pressure, and $(L_t)_{t\geq 0}$
is an infinite dimensional cylindrical L\'evy process given by
$$
L_t=\sum_{j\in\mN}\beta_j L^{(j)}_te_j,
$$
where $\{(L^{(j)}_t)_{t\geq 0},j\in\mN\}$ is a sequence of independent one dimensional
purely discontinuous L\'evy processes defined on some filtered probability space
$(\Omega,\sF,(\sF_t)_{t\geq 0}; P)$ and with the same L\'evy measure $\nu$,
$\{\beta_j,j\in\mN\}$ is a sequence of positive numbers and $\{e_j,j\in\mN\}$
is a sequence of orthogonal basis of Hilbert space $\mH^0$, where for $\gamma\in\mR$,
$\mH^\gamma$ with the norm $\|\cdot\|_\gamma$ and inner product $\<\cdot,\cdot\>_\gamma$
denotes the usual Sobolev space of divergence free vector fields on $\mT^2$ (see Section 2 for a definition).

As a continuous model, stochastic Navier-Stokes equation driven by
Brownian motion has been extensively studied in the past decades (cf. \cite{Ha-Ma,Da-De,De-Od,Fl-Ro}, etc.).
Meanwhile,  stochastic partial differential equation with jump has also been studied recently
(cf. \cite{Pe-Za,Do-Zh}). However, in the well-known results, the assumption that the jump process has finite second order moments
was required in order to obtain the usual energy estimate. This excludes the interest $\alpha$-stable process.
In this note, we establish the well posedness for stochastic 2D Navier-Stokes equation (\ref{NSE})
driven by a general cylindrical L\'evy process,
and obtain the existence of invariant measures for this discontinuous model. More precisely, we shall prove that:
\bt\label{Main}
Suppose that for some $\theta\in(0,1]$,
$$
\mbox{\bf (H$_\theta$):}\quad H_\theta:=\int_{|x|>1}|x|^\theta\nu(\dif x)+\sum_{j\in\mN}|\beta_j|^\theta<+\infty.
$$
Then for any $\varphi\in\mH^0$,
there exists a unique solution $(u_t)_{t\geq 0}=(u_t(\varphi))_{t\geq 0}$
to equation(\ref{NSE}) satisfying that for almost all $\omega$ and  for any $t>0$,
\begin{enumerate}[(i)]
\item $t\mapsto u_t(\omega)$ is right continuous and has left-hand limit in $\mH^0$, and
$\int^t_0\|\nabla u_s(\omega)\|^2_0\dif s<+\infty$;
\item  it holds that for any $\phi\in\mH^1$,
$$
\<u_t(\omega),\phi\>_0=\<\varphi,\phi\>_0+\int^t_0[\<\Delta u_s(\omega),\phi\>_0
+\<u_s(\omega)\otimes u_s(\omega),\nabla\phi\>_0]\dif s+\<L_t(\omega),\phi\>_0.
$$
\end{enumerate}
Moreover, there exists a constant $C=C(H_\theta,\theta)>0$ such that for any $t>0$,
\begin{align}
\mE\left(\sup_{s\in[0,t]}\|u_s\|^\theta_0\right)+
\mE\left(\int^t_0\frac{\|\nabla u_s\|^2_0}{(\|u_s\|^2_0+1)^{1-\theta/2}}\dif s\right)
\leq C(1+\|\varphi\|_0^\theta+t).\label{PL9}
\end{align}
In particular, there exists a probability measure $\mu$ on $(\mH^0,\sB(\mH^0))$
called invariant measure of $(u_t(\varphi))_{t\geq 0}$ such that
for any bounded measurable function $\Phi$ on $\mH^0$,
$$
\int_{\mH^0}\mE \Phi(u_t(\varphi))\mu(\dif\varphi)=\int_{\mH^0} \Phi(\varphi)\mu(\dif\varphi).
$$
\et
\br
Assumption {\bf (H$_\theta$)} implies that cylindrical L\'evy process $(L_t)_{t\geq 0}$ admits a cadlag version in $\mH^0$ and
for any $t>0$ (cf. \cite[p.159, Theorem 25.3]{Sa}),
$$
\mE\|L_t\|^\theta_0<+\infty.
$$
In fact, for $\theta\in(0,1]$,
by the elementary inequality $(a+b)^\theta\leq a^\theta+b^\theta$, we have
$$
\mE\|L_t\|^\theta_0
\leq\mE\left(\sum_{j\in\mN}|\beta_j|\cdot|L^{(j)}_t|\right)^\theta
\leq \sum_{j\in\mN}|\beta_j|^\theta\cdot\mE|L^{(j)}_t|^\theta
=\mE|L^{(1)}_t|^\theta\sum_{j\in\mN}|\beta_j|^\theta<+\infty.
$$
\er
\br\label{Re1}
By estimate (\ref{PL9}) and Poinc\`are's inequality, we have
\begin{align*}
\mE\left(\int^t_0\|\nabla u_s\|^\theta_0\dif s\right)&\leq
\mE\left(\int^t_0\frac{\|\nabla u_s\|^\theta_0(\|u_s\|^{2-\theta}_0+1)}{(\|u_s\|^2_0+1)^{1-\theta/2}}\dif s\right)\\
&\leq C\mE\left(\int^t_0\frac{\|\nabla u_s\|^2_0+1}{(\|u_s\|^2_0+1)^{1-\theta/2}}\dif s\right)\\
&\leq C(1+\|\varphi\|_0^\theta+t).
\end{align*}
This estimate in particular yields the existence of invariant
measures by the classical Bogoliubov-Krylov's argument (cf.
\cite{Da-Za}). \er \br An obvious open question is about the
uniqueness of invariant measures (i.e. ergodicity) for discontinuous
system (\ref{NSE}). The notion of asymptotic strong Feller property
in \cite{Ha-Ma} is perhaps helpful for solving this problem. \er
This paper is organized as follows: In Section 2, we give some
necessary materials. In Section 3, we prove the main result. 
\section{Preliminaries}

In this section we prepare some materials for later use. Let $C^\infty_0(\mT^2)^2$ be the space of all
smooth $\mR^2$-valued function on $\mT^2$ with vanishing mean and divergence, i.e.,
$$
\int_{\mT^2}f(x)\dif x=0,\ \ \div f(x)=0.
$$
For $\gamma\in\mR$, let $\mH^\gamma$ be the completion of $C^\infty_0(\mT^2)^2$ with respect to the norm
$$
\|f\|_\gamma=\left(\int_{\mT^2}|(-\Delta)^{\gamma/2}f(x)|^2\dif x\right)^{1/2},
$$
where $(-\Delta)^{\gamma/2}$ is defined through Fourier's transform. Thus, $(\mH^\gamma,\|\cdot\|_\gamma)$ is
a separable Hilbert space with the obvious inner product
$$
\<f,g\>_\gamma:=\int_{\mT^2}(-\Delta)^{\gamma/2}f(x)\cdot(-\Delta)^{\gamma/2}g(x)\dif x.
$$
Below, we shall fix an orthogonal basis $\{e_j, j\in\mN\}\subset C^\infty_0(\mT^2)^2$ of $\mH^0$ consisting of the eigenvectors of $\Delta$, i.e.,
\begin{align}
\Delta e_j=-\lambda_j e_j, \ \ \<e_j,e_j\>_0=1, \ \ j=1,2,\cdots,
\label{PL7}\end{align}
where $0<\lambda_1<\cdots<\lambda_j\uparrow\infty$.

Let $\{(L^{(j)}_t)_{t\geq 0},j\in\mN\}$ be a sequence of independent one dimensional purely discontinuous L\'evy processes
with the same characteristic function, i.e.,
$$
\mE e^{\mathrm{i}\xi L^{(j)}_t}=e^{-t\psi(\xi)},\ \forall t\geq 0, j=1,2,\cdots,
$$
where $\psi(\xi)$ is a complex valued function called L\'evy symbol given by
$$
\psi(\xi)=\int_{\mR\setminus\{0\}}(e^{\mathrm{i}\xi y}-1-\mathrm{i}\xi y1_{|y|\leq 1})\nu(\dif y),
$$
where $\nu$ is the L\'evy measure and satisfies that
$$
\int_{\mR\setminus\{0\}}1\wedge |y|^2\nu(\dif y)<+\infty.
$$
For $t>0$ and $\Gamma\in\sB(\mR\setminus\{0\})$,
the Poisson random measure associated with $L^{(j)}_t$ is defined by
$$
N^{(j)}(t,\Gamma):=\sum_{s\in(0,t]}1_{\Gamma}(L^{(j)}_s-L^{(j)}_{s-}).
$$
The compensated Poisson random measure is given by
$$
\tilde N^{(j)}(t,\Gamma)=N^{(j)}(t,\Gamma)-t\nu(\Gamma).
$$
By L\'evy-It\^o's decomposition (cf. \cite[p.108, Theorem 2.4.16]{Ap}), one has
$$
L^{(j)}_t=\int_{|x|\leq1}x\tilde N^{(j)}(t,\dif x)+\int_{|x|>1}x N^{(j)}(t,\dif x).
$$

For a Polish space $(\mG,\rho)$, let $\mD(\mR_+;\mG)$ be the space of all
right continuous functions with left-hand limits from $\mR_+$
to $\mG$, which is endowed with the Skorohod topology:
\begin{align}
d_\mG(u,v):=\inf_{\lambda\in\Lambda}\left[\sup_{s\not= t}\left|\log\frac{\lambda(t)-\lambda(s)}{t-s}\right|
\vee\int^\infty_0 \sup_{t\geq 0}(\rho(u_{t\wedge r},v_{\lambda(t)\wedge r})\wedge 1)e^{-r}\dif r\right],
\label{PL6}
\end{align}
where $\Lambda$ is the space of all continuous and strictly increasing function from $\mR_+\to\mR_+$
with $\lambda(0)=0$ and $\lambda(\infty)=\infty$. Thus,
$(\mD(\mR_+;\mG);d)$ is again a Polish space (cf. \cite[p.121, Theorem 5.6]{Ei-Ku}).

We need the following tightness criterion,
which is a direct combination of \cite[Corollary 5.2]{Ja} and Aldous's criterion \cite{Al}.

\bt\label{Th1}
Let $\{(X^n_t)_{t\geq 0}, n\in\mN\}$ be a sequence of $\mH^{-1}$-valued stochastic processes with cadlag path. Assume that
\begin{enumerate}[(i)]
\item for each $\phi\in C^\infty_0(\mT^2)^2$ and $t>0$,
$\lim_{K\to\infty}\sup_{n\in\mN}P\Big\{\sup_{s\in[0,t]}|\<X^n_s, \phi\>_{-1}|\geq K\Big\}=0$;

\item for each $\phi\in C^\infty_0(\mT^2)^2$ and $t,a>0$, $\lim_{\eps\to 0+}\sup_{n\in\mN}\sup_{\tau\in\cS_t}
P\Big\{|\<X^n_\tau-X^n_{\tau+\eps}, \phi\>_{-1}|\geq a\Big\}=0$, where $\cS_t$ denotes all the
($\sF_t$)-stopping times with bound $t$;
\item for every $\eps>0$ and $t>0$,
$$
\lim_{m\to\infty}\sup_{n\in\mN} P\left(\sup_{s\in[0,t]}\sum_{j=m}^\infty\<X^n_s,e_j\>^2_{-1}\geq\eps\right)=0.
$$
\end{enumerate}
Then the laws of $(X^n_t)_{t\geq 0}$ in $\mD(\mR_+;\mH^{-1})$ is tight.
\et

The following result comes from \cite[p.131 Theorem 7.8]{Ei-Ku}.
\bt\label{Th2}
Suppose that stochastic processes sequence $\{(X^n_t)_{t\geq 0}, n\in\mN\}$ weakly converges to
$(X_t)_{t\geq 0}$  in $\mD(\mR_+;\mH^{-1})$.Then,
for any $t>0$ and $\phi\in\mH^1$, there exists a sequence $t_n\downarrow t$ such that
for any bounded continuous function $f$,
$$
\lim_{n\to\infty}\mE f(\<X^n_{t_n},\phi\>_{-1})=\mE f(\<X_t,\phi\>_{-1}).
$$
\et

We also need the following technical result.
\bl\label{Le1}
Suppose that sequence $\{u^n,n\in\mN\}$ converges to $u$ in $\mD(\mR_+;\mH^{-1})$. Then for any $T>0$ and $m\in\mN$,
\begin{align}
\sup_{t\in[0,T]}\|u_t\|_0\leq \varliminf_{n\to\infty}\sup_{t\in[0,T+\frac{1}{m}]}\|u^n_t\|_0.\label{PL5}
\end{align}
If in addition, for Lebesgue almost all $t$, $u^n_t$ converges to $u_t$ in $\mH^0$, then for any $\beta>0$,
\begin{align}
\int^T_0\frac{\|\nabla u_t\|^2_0}{(1+\|u_t\|^2_0)^\beta}\dif t\leq \varliminf_{n\to\infty}
\int^T_0\frac{\|\nabla u^n_t\|^2_0}{(1+\|u^n_t\|^2_0)^\beta}\dif t.\label{PL55}
\end{align}
\el
\begin{proof}
Without loss of generality, we assume that the right hand side of (\ref{PL5}) is finite.
For any $\phi\in\mH^1$, it is clear that  $t\mapsto\<u_t,\phi\>_0$ is a cadlag real valued function,
and by definition (\ref{PL6}) of Skorohod metric, we have
$$
d_{\mR}(\<u^n,\phi\>_0, \<u,\phi\>_0)\leq (2+\|\phi\|_1)d_{\mH^{-1}}(u^n,u),
$$
and so $\<u^n,\phi\>_0$ converges to $\<u,\phi\>_0$ in $\mD(\mR_+;\mR)$ as $n\to\infty$.
Since the discontinuous points of $\<u_\cdot,\phi\>_0$ are at most countable, for any $T>0$ and $m\in\mN$, there
exists a time $T_m\in(T,T+1/m)$ such that $\<u_\cdot,\phi\>_0$ is continuous at $T_m$. Thus,
we have (cf. \cite[p.119, Proposition 5.3]{Ei-Ku})
$$
\lim_{n\to\infty}\sup_{t\in[0,T_m]}|\<u^n_t,\phi\>_0|=\sup_{t\in[0,T_m]}|\<u_t,\phi\>_0|.
$$
Hence,
\begin{align*}
\sup_{t\in[0,T]}\|u_t\|_0&=\sup_{t\in[0,T]}\sup_{\phi\in\mH^1;\|\phi\|_0\leq 1}|\<u_t,\phi\>_0|\\
&\leq\sup_{\phi\in\mH^1;\|\phi\|_0\leq 1}\sup_{t\in[0,T_m]}|\<u_t,\phi\>_0|\\
&=\sup_{\phi\in\mH^1;\|\phi\|_0\leq 1}\lim_{n\to\infty}\sup_{t\in[0,T_m]}|\<u^n_t,\phi\>_0|\\
&\leq\varliminf_{n\to\infty}\sup_{\phi\in\mH^1;\|\phi\|_0\leq 1}\sup_{t\in[0,T_m]}|\<u^n_t,\phi\>_0|\\
&=\varliminf_{n\to\infty} \sup_{t\in[0,T_m]}\|u^n_t\|_0.
\end{align*}
Thus, (\ref{PL5}) is proven.

For proving (\ref{PL55}), let $\cN$ be the Lebesgue null set such that for all $t\notin \cN$,
$u^n_t$ converges to $u_t$ in $\mH^0$. Fixing a $t\notin\cN$, then as above, we have
$$
\frac{\|\nabla u_t\|^2_0}{(1+\|u_t\|^2_0)^\beta}\leq
\frac{\varliminf_{n\to\infty}\|\nabla u^n_t\|^2_0}{(1+\lim_{n\to\infty}\|u^n_t\|^2_0)^\beta}\leq
\varliminf_{n\to\infty}\frac{\|\nabla u^n_t\|^2_0}{(1+\|u^n_t\|^2_0)^\beta}.
$$
Estimate (\ref{PL55}) now follows by Fatou's lemma.
\end{proof}

\section{Proof of Theorem \ref{Main}}

We first give the following definition about the weak solutions to equation (\ref{NSE}).
\bd\label{Def1}
A probability measure $P$ on $\mD(\mR_+;\mH^{-1})$ is called a weak solution of equation (\ref{NSE}) if
\begin{enumerate}[(i)]
\item for any $t>0$,
$P\left(u\in\mD(\mR_+;\mH^{-1}): \sup_{s\in[0,t]}\|u_s\|_0+\int^t_0\|\nabla u_s\|^2_0\dif s<+\infty\right)=1$;
\item for any $j\in\mN$,
\begin{align}
M^{(j)}_t(u):=\<u_t,e_j\>_0-\<u_0,e_j\>_0-
\int^t_0[\<u_s,\Delta e_j\>_0+\<u_s\otimes u_s,\nabla e_j\>_0]\dif s\label{Levy}
\end{align}
is a L\'evy process with the characteristic function
$$
\mE e^{\mathrm{i}\xi M^{(j)}_t}
=\exp\left\{t\int_{\mR\setminus\{0\}}(e^{\mathrm{i}\xi y\beta_j}-1
-\mathrm{i}\xi y \beta_j1_{|y|\leq 1})\nu(\dif y)\right\},
$$
and $\{(M^{(j)}_t)_{t\geq 0}, j\in\mN\}$ is a sequence of independent L\'evy processes.
\end{enumerate}
\ed

{\it Proof of Existence of Weak Solutions}: We use Galerkin's approximation to prove
the existence of weak solutions and divide the proof into three steps.

(Step 1): For $n\in\mN$, set
$$
\mH^0_n:=\mathrm{span}\{e_1,e_2,\cdots,e_n\},
$$
and let $\Pi_n$ be the projection from $\mH^0$ to $\mH^0_n$ and define
$$
L^n_t:=\sum_{j=1}^n\beta_jL^{(j)}_te_j=\sum_{j=1}^n\int_{|y|\leq1}y\beta_je_j\tilde N^{(j)}(t,\dif y)
+\sum_{j=1}^n\int_{|y|>1}y\beta_je_j N^{(j)}(t,\dif y).
$$
Consider the following finite dimensional SDE driven by finite dimensional L\'evy process $L^n_t$:
\begin{align}
\dif u^n_t=[\Delta u^n_t-\Pi_n((u^n_t\cdot\nabla) u^n_t)]\dif t+\dif L^n_t,\ \ u^n_0=\Pi_n\varphi.\label{SDE}
\end{align}
Since for any $R>0$ and $u,v\in\mH^0_n$ with $\|u\|_0,\|v\|_0\leq R$,
$$
\|\Pi_n((u\cdot\nabla)u-(v\cdot\nabla)v)\|_0\leq C_{R,n}\|u-v\|_0
$$
and
\begin{align}
\<u,\Delta u-\Pi_n((u\cdot\nabla)u)\>_0=-\|\nabla u\|_0, \ \ \forall u\in\mH^0_n,\label{E2}
\end{align}
finite dimensional SDE (\ref{SDE}) is thus well-posed.

Define a smooth function $f_n$ on $\mH^0_n$ by
$$
f_n(u):=(\|u\|^2_0+1)^{\theta/2},\ \ u\in\mH^0_n.
$$
By simple calculations, we have
\begin{align}
\nabla f_n(u)=\frac{\theta u}{(\|u\|^2_0+1)^{1-\theta/2}},\ \
\nabla^2 f_n(u)=\frac{\theta\sum_{i=1}^ne_i\otimes e_i}{(\|u\|^2_0+1)^{1-\theta/2}}
-\frac{\theta(2-\theta)u\otimes u}{(\|u\|^2_0+1)^{2-\theta/2}},\label{E1}
\end{align}
and for all $u,v\in\mH^0_n$,
\begin{align}
|f_n(u)-f_n(v)|\leq|(\|u\|^2_0+1)^{1/2}-(\|v\|^2_0+1)^{1/2}|^\theta\leq \|u-v\|^\theta_0.\label{E4}
\end{align}
By (\ref{SDE}), (\ref{E2}), (\ref{E1}) and It\^o's formula (cf. \cite[p.226, Theorem 4.4.7]{Ap}), we have
\begin{align*}
f_n(u^n_t)&=f_n(u^n_0)-\int^t_0\frac{\theta\|\nabla u^n_s\|^2_0}{(\|u^n_s\|^2_0+1)^{1-\theta/2}}\dif s
+\sum_{j=1}^n\int^t_0\!\!\!\int_{|y|\leq 1}[f_n(u^n_s+y\beta_j e_j)-f_n(u^n_s)]\tilde N^{(j)}(\dif s,\dif y)\\
&\quad+\sum_{j=1}^n\int^t_0\!\!\!\int_{|y|\leq 1}\left[f_n(u^n_s+y\beta_j e_j)-f_n(u^n_s)
-\frac{\theta\<u^n_s,y\beta_je_j\>_0}{(|u^n_s|^2+1)^{1-\theta/2}}\right]\nu(\dif y)\dif s\\
&\quad+\sum_{j=1}^n\int^t_0\!\!\!\int_{|y|> 1}\left[f_n(u^n_s+y\beta_j e_j)-f_n(u^n_s)\right]N^{(j)}(\dif s,\dif y)\\
&=:f_n(u^n_0)-I^n_1(t)+I^n_2(t)+I^n_3(t)+I^n_4(t).
\end{align*}
For $I^n_2(t)$, by Burkholder's inequality and (\ref{E4}), we have
\begin{align*}
\mE\left(\sup_{t\in[0,T]}I^n_2(t)\right)&\leq C\sum_{j=1}^n\mE\left(\int^T_0\!\!\!\int_{|y|\leq 1}
|f_n(u^n_s+y\beta_j e_j)-f_n(u^n_s)|^2N^{(j)}(\dif s,\dif y)\right)^{1/2}\\
&\leq C\sum_{j=1}^n\left(\mE\int^T_0\!\!\!\int_{|y|\leq 1}
|f_n(u^n_s+y\beta_j e_j)-f_n(u^n_s)|^2\nu(\dif y)\dif s\right)^{1/2}\\
&\leq C T^{1/2}\sum_{j=1}^n|\beta_j|^\theta\left(\int_{|y|\leq 1}|y|^{2\theta}\nu(\dif y)\right)^{1/2}\\
&\leq C T^{1/2}\sum_{j=1}^\infty|\beta_j|^\theta\left(\int_{|y|\leq 1}|y|^2\nu(\dif y)\right)^{1/2}.
\end{align*}
where we have used condition {\bf (H$_\theta$)}. Here and after,
the constant $C$ is independent of $n$ and $T$.
For $I^n_3(t)$, by Taylor's expansion and (\ref{E1}), we have
\begin{align*}
\mE\left(\sup_{t\in[0,T]}I^n_3(t)\right)
\leq C\sum_{j=1}^n\beta_j^2\int^T_0\!\!\!\int_{|y|\leq 1}|y|^2\nu(\dif y)\dif s
\leq CT\sum_{j=1}^\infty|\beta_j|^\theta \int_{|y|\leq 1}|y|^2\nu(\dif y).
\end{align*}
For $I^n_4(t)$, by (\ref{E4}), we have
\begin{align*}
\mE\left(\sup_{t\in[0,T]}I^n_4(t)\right)&\leq \sum_{j=1}^n\mE\left(\int^T_0\!\!\!\int_{|y|> 1}|f_n(u^n_s+y\beta_je_j)-f_n(u^n_s)|N^{(j)}(\dif s,\dif y)\right)\\
&=\sum_{j=1}^n\mE\left(\int^T_0\!\!\!\int_{|y|> 1}|f_n(u^n_s+y\beta_j e_j)-f_n(u^n_s)|\nu(\dif y)\dif s\right)\\
&\leq C T\sum_{j=1}^\infty|\beta_j|^\theta\int_{|y|> 1}|y|^\theta\nu(\dif y).
\end{align*}
Combining the above calculations, we obtain that
\begin{align}
\mE\left(\sup_{t\in[0,T]}(\|u^n_t\|^2_0+1)^{\theta/2}\right)
+\mE\int^T_0\frac{\theta\|\nabla u^n_s\|^2_0}{(\|u^n_s\|^2_0+1)^{1-\theta/2}}\dif s\leq
(\|\varphi\|^2_0+1)^{\theta/2}+CT+CT^{1/2}.\label{PL3}
\end{align}

(Step 2): In this step, we use Theorem \ref{Th1} to show
that $\{(u^n_t)_{t\geq 0},n\in\mN\}$ is tight in $\mD(\mR_+;\mH^{-1})$.
For any $\phi\in C^\infty_0(\mT^2)^2$, by equation (\ref{SDE}), we have
\begin{align*}
\<u^n_t,\phi\>_{-1}&=\<u^n_0,\phi\>_{-1}+\int^t_0[\<\Delta u^n_s,\phi\>_{-1}-\<(u^n_s\cdot\nabla) u^n_s,\phi\>_{-1}]\dif s
+\<L^n_t,\phi\>_{-1}\\
&=\<u^n_0,\phi\>_{-1}+\int^t_0[\<u^n_s,\Delta \phi\>_{-1}+\<u^n_s\otimes u^n_s,\nabla\phi\>_{-1}]\dif s
+\<L^n_t,\phi\>_{-1}.
\end{align*}
Thus, for $\eps>0$ and any stopping time $\tau$ bounded by $t$, we have
\begin{align*}
\<u^n_{\tau+\eps}-u^n_\tau,\phi\>_{-1}&=\int^{\tau+\eps}_\tau
[\<u^n_s,\Delta \phi\>_{-1}+\<u^n_s\otimes u^n_s,\nabla\phi\>_{-1}]\dif s
+\<L^n_{\tau+\eps}-L^n_\tau,\phi\>_{-1}\\
&\leq \eps\sup_{s\in[0,t]}\left(\|u^n_s\|_0\cdot\|\phi\|_0+\|u^n_s\|^2_0\cdot\|\nabla(-\Delta)^{-1}\phi\|_\infty\right)\\
&\quad+\sum_{j=1}^n|\beta_j|\cdot|L^{(j)}_{\tau+\eps}-L^{(j)}_\tau|\cdot\|(-\Delta)^{-1}\phi\|_0.
\end{align*}
Using $(a+b)^\theta\leq a^\theta+b^\theta$ provided that $\theta\in(0,1]$, we get
\begin{align*}
\mE|\<u^n_{\tau+\eps}-u^n_\tau,\phi\>_{-1}|^{\theta/2}&\leq C_\phi\mE\left(\sup_{s\in[0,t]}\|u^n_s\|^\theta_0+1\right)\eps^{\theta/2}
+C_\phi\left(\mE\sum_{j=1}^n|\beta_j|^\theta\cdot|L^{(j)}_{\tau+\eps}-L^{(j)}_\tau|^\theta\right)^{1/2}
\end{align*}
By the strong Markov property of L\'evy process (cf. \cite[p.278, Theorem 40.10]{Sa}), we have
$$
\mE|L^{(j)}_{\tau+\eps}-L^{(j)}_\tau|^\theta=\mE|L^{(j)}_\eps|^\theta=\mE|L^{(1)}_\eps|^\theta,\ \forall j\in\mN.
$$
Thus, by (\ref{PL3}) and {\bf (H$_\theta$)},
\begin{align}
\mE|\<u^n_{\tau+\eps}-u^n_\tau,\phi\>_{-1}|^{\theta/2}\leq C\Big[\eps^{\theta/2}
+(\mE|L^{(1)}_\eps|^\theta)^{1/2}\Big],\label{PL1}
\end{align}
where the constant $C$ is independent of $n,\tau$ and $\eps$.
On the other hand, by (\ref{PL7}), we have
\begin{align}
\mE\left(\sup_{s\in[0,t]}\sum_{j=m}^\infty\<u^n_s,e_j\>_{-1}^2\right)^{\theta/2}
=\mE\left(\sup_{s\in[0,t]}\sum_{j=m}^\infty\frac{\<u^n_s,e_j\>_0^2}{\lambda^2_j}\right)^{\theta/2}
\leq\frac{1}{\lambda^\theta_m}\mE\left(\sup_{s\in[0,t]}\|u^n_s\|^\theta_0\right).\label{PL2}
\end{align}
By Theorem \ref{Th1} and (\ref{PL3})-(\ref{PL2}),
one knows that the law of $(u^n_t)_{t\geq 0}$ in
$\mD(\mR_+;\mH^{-1})$ denoted by $P_n$ is tight.

(Step 3): Let $P$ be any accumulation point of $\{P_n,n\in\mN\}$.
In this step, we show that $P$ is a weak solution of equation (\ref{NSE}) in the sense of Definition
\ref{Def1}. First of all, by Skorohod's embedding theorem,
there exists a probability space $(\tilde\Omega,\tilde\sF,\tilde P)$ and $\mD(\mR_+;\mH^{-1})$-valued
random variables $X^n$ and $X$ such that

(i) Law of $X^n$ under $\tilde P$ is $P_n$ and law of $X$ under $\tilde P$ is $P$.

(ii) $X^n$ converges to $X$ in $\mD(\mR_+;\mH^{-1})$ a.s. as $n\to\infty$.

Thus, by (\ref{PL3}), we have
\begin{align}
\tilde\mE\left(\sup_{t\in[0,T]}\|X^n_t\|_0^\theta\right)
+\tilde\mE\left(\int^T_0\frac{\theta\|\nabla X^n_s\|^2_0}{(\|X^n_s\|^2_0+1)^{1-\theta/2}}\dif s\right)\leq
C(1+\|\varphi\|^\theta_0+T).\label{PL30}
\end{align}
By Lemma \ref{Le1} and Fatou's lemma, for any $m\in\mN$, we have
\begin{align}
\mE^P\left(\sup_{t\in[0,T]}\|u_t\|^\theta_0\right)
&=\tilde\mE\left(\sup_{t\in[0,T]}\|X_t\|^\theta_0\right)
\leq\varliminf_{n\to\infty}\tilde\mE\left(\sup_{t\in[0,T+1/m]}\|X^n_t\|^\theta_0\right)\no\\
&\leq (\|\varphi\|^2_0+1)^{\theta/2}+C(T+1/m)+C(T+1/m)^{1/2}.\label{Op1}
\end{align}
On the other hand, for any $\delta\in(0,\theta/4)$, by H\"older's inequality and (\ref{PL30}), we have
\begin{align*}
\tilde\mE\left(\int^T_0\|X^n_s-X_s\|_0^\delta\dif s\right)
&\leq\tilde\mE\left(\int^T_0\|X^n_s-X_s\|^{\delta/2}_{-1}\|X^n_s-X_s\|^{\delta/2}_1\dif s\right)\\
&\leq\left(\tilde\mE\int^T_0\|X^n_s-X_s\|^{\delta}_{-1}\dif s\right)^{1/2}
\left(\tilde\mE\int^T_0\|X^n_s-X_s\|^{\delta}_1\dif s\right)^{1/2}\to 0.
\end{align*}
So, there exists a subsequence still denoted by $n$ such that for $\tilde P\times\dif t$-almost all $(\omega,s)$,
$X^n_s(\omega)$ converges to $X_s(\omega)$ in $\mH^0$. By Lemma \ref{Le1} and (\ref{PL30}), we then obtain
\begin{align}
\mE^P\left(\int^T_0\frac{\theta\|\nabla u_s\|^2_0}{(\|u_s\|^2_0+1)^{1-\theta/2}}\dif s\right)
&=\tilde\mE\left(\int^T_0\frac{\theta\|\nabla X_s\|^2_0}{(\|X_s\|^2_0+1)^{1-\theta/2}}\dif s\right)\no\\
&\leq\varliminf_{n\to\infty}\tilde\mE\left(\int^T_0\frac{\theta\|\nabla X^n_s\|^2_0}{(\|X^n_s\|^2_0+1)^{1-\theta/2}}\dif s\right)\no\\
&\leq C(1+\|\varphi\|^\theta_0+T).\label{Op2}
\end{align}
Combining (\ref{Op1}) and (\ref{Op2}) gives (\ref{PL9}). In particular,
$\sup_{t\in[0,T]}\|u_t\|_0$ and $\int^T_0\frac{\theta\|\nabla u_s\|^2_0}{(\|u_s\|^2_0+1)^{1-\theta/2}}\dif s$
is finite $P$-almost surely, which produces (i) of Definition \ref{Def1}.

Fixing $j\in\mN$, in order to show that $M^{(j)}_t$ defined by (\ref{Levy}) is
a L\'evy process, we only need to prove that for any $0\leq s<t$,
\begin{align}
\mE^P e^{\mathrm{i} \xi(M^{(j)}_t-M^{(j)}_s)}=\tilde \mE e^{\mathrm{i} \xi(\tilde M^{(j)}_t-\tilde M^{(j)}_s)}
=\exp\left\{(t-s)\int_{\mR\setminus\{0\}}(e^{\mathrm{i}\xi y\beta_j}-1
-1_{|y|\leq 1}\mathrm{i}\xi y\beta_j)\nu(\dif y)\right\},\label{PL8}
\end{align}
where
$$
\tilde M^{(j)}_t:=\<X_t,e_j\>_0-\<X_0,e_j\>_0-
\int^t_0[\<X_r,\Delta e_j\>_0+\<X_r\otimes X_r,\nabla e_j\>_0]\dif r.
$$
Fix $0\leq s<t$ below. By Theorem \ref{Th2}, there exists $(s_n,t_n)\downarrow (s,t)$ such that
$$
\lim_{n\to\infty}\tilde\mE e^{\mathrm{i} \xi \<X^n_{t_n},e_j\>_0}=\tilde \mE e^{\mathrm{i} \xi \<X_t,e_j\>_0},\ \
\lim_{n\to\infty}\tilde \mE e^{\mathrm{i} \xi \<X^n_{s_n},e_j\>_0}=\tilde\mE e^{\mathrm{i} \xi \<X_s,e_j\>_0}.
$$
By equation (\ref{SDE}), it is well-known that for any $n\geq j$,
\begin{align*}
&\tilde\mE \exp\left\{\mathrm{i} \xi \left[\<X^n_{t_n}-X^n_{s_n},e_j\>_0
-\int^{t_n}_{s_n}[\<X^n_r,\Delta e_j\>_0+\<X^n_r\otimes X^n_r,\nabla e_j\>_0]\dif r\right]\right\}\\
&\quad=\mE^{P_n} \exp\left\{\mathrm{i} \xi\left[\<u^n_{t_n}-u^n_{s_n},e_j\>_0
-\int^{t_n}_{s_n}[\<u^n_r,\Delta e_j\>_0
+\<u^n_r\otimes u^n_r,\nabla e_j\>_0]\dif r\right]\right\}\\
&\quad=\exp\left\{(t_n-s_n)\int_{\mR\setminus\{0\}}(e^{\mathrm{i}\xi y\beta_j}-1
-1_{|y|\leq 1}\mathrm{i}\xi y\beta_j)\nu(\dif y)\right\}.
\end{align*}
Thus, for proving (\ref{PL8}), it suffices to prove the following limits:
\begin{align*}
&\lim_{n\to\infty}\tilde\mE\left|\exp\left\{\mathrm{i}\xi\int^t_s\<X^n_r\otimes X^n_r,\nabla e_j\>_0\dif r\right\}
-\exp\left\{\mathrm{i}\xi\int^t_s\<X_r\otimes X_r,\nabla e_j\>_0\dif r\right\}\right|=0,\\
&\lim_{n\to\infty}\tilde\mE\left|\exp\left\{\mathrm{i}\xi\int^t_s\<X^n_r,\Delta e_j\>_0\dif r\right\}
-\exp\left\{\mathrm{i}\xi\int^t_s\<X_r,\Delta e_j\>_0\dif r\right\}\right|=0,\\
&\lim_{n\to\infty}\tilde\mE\left|\exp\left\{\mathrm{i}\xi\int^{t_n}_{s_n}\<X^n_r\otimes X^n_r,\nabla e_j\>_0\dif r\right\}
-\exp\left\{\mathrm{i}\xi\int^t_s\<X^n_r\otimes X^n_r,\nabla e_j\>_0\dif r\right\}\right|=0,\\
&\lim_{n\to\infty}\tilde\mE\left|\exp\left\{\mathrm{i}\xi\int^{t_n}_{s_n}\<X^n_r,\Delta e_j\>_0\dif r\right\}
-\exp\left\{\mathrm{i}\xi\int^t_s\<X^n_r,\Delta e_j\>_0\dif r\right\}\right|=0.
\end{align*}
Let us only prove the first limit, the others are similar. Noticing that for any $\delta\in(0,1)$ and $a,b\in\mR$,
$$
|e^{\mathrm{i}a}-e^{\mathrm{i}b}|\leq 2(|a-b|\wedge 1)\leq 2|a-b|^\delta,
$$
by H\"older's inequality and $\|u\|_0\leq\|u\|_{-1}^{1/2}\|u\|_1^{1/2}$, we have for $\delta<\theta/4$,
\begin{align*}
&\tilde\mE\left|\exp\left\{\mathrm{i}\xi\int^t_s\<X^n_r\otimes X^n_r,\nabla e_j\>_0\dif r\right\}
-\exp\left\{\mathrm{i}\xi\int^t_s\<X_r\otimes X_r,\nabla e_j\>_0\dif r\right\}\right|\\
&\quad\leq2|\xi|^\delta\tilde\mE\left|\int^t_s\<X^n_r\otimes X^n_r-X_r\otimes X_r,
\nabla e_j\>_0\dif r\right|^\delta\\
&\quad\leq C \tilde\mE\left(\int^t_s\|X^n_r-X_r\|_0(\|X^n_r\|_0+\|X_r\|_0)\dif r\right)^\delta\\
&\quad\leq C \tilde\mE\left(\sup_{r\in[s,t]}(\|X^n_r\|_0+\|X_r\|_0)
\int^t_s\|X^n_r-X_r\|^{1/2}_{-1}\|X^n_r-X_r\|^{1/2}_1\dif r\right)^\delta\\
&\quad\leq C \tilde\mE\Bigg(\sup_{r\in[s,t]}(\|X^n_r\|_0+\|X_r\|_0+1)^{2\delta-(\theta\delta/2)}
\left(\int^t_s\|X^n_r-X_r\|_{-1}\dif r\right)^{\delta/2}\\
&\qquad\times\left(\int^t_s\frac{(\|X^n_r\|_1+\|X_r\|_1)}{(\|X^n_r\|^2_0+\|X_r\|^2_0+1)^{1-\theta/2}}
\dif r\right)^{\delta/2}\Bigg)\\
&\quad\leq C\left[\tilde\mE\left(\int^t_s\|X^n_r-X_r\|_{-1}\dif r\right)^{2\delta}\right]^{1/4}\to 0,
\end{align*}
as $n\to\infty$, where in the last inequality, we have used (\ref{PL30}) and H\"older's inequality.
As for the independence of $M^{(j)}$ for different $j\in\mN$, it can be proved in a similar way.

\vspace{5mm}

{\it Proof of Theorem \ref{Main}}: The pathwise uniqueness follows  by the classical result for 2D
deterministic Navier-Stokes equation.
As for the existence of invariant measures, basing on (\ref{PL9}) (see Remark \ref{Re1}),
it follows by the classical Bogoliubov-Krylov's argument.

\end{document}